\documentclass{article}
\usepackage[utf8]{inputenc}
\usepackage{authblk}
\usepackage{setspace}
\usepackage{geometry}
\usepackage{enumitem}
\usepackage[table]{xcolor}
\usepackage{tabularx}
\usepackage{graphicx}
\graphicspath{ {./figures/} }
\usepackage{subcaption}
\usepackage{amsmath}
\usepackage{lineno}
\usepackage{amssymb}
\usepackage{amsthm}
\usepackage{algorithm}
\usepackage{comment}
\makeatletter
\usepackage{booktabs}

\usepackage[noend]{algpseudocode}
\usepackage{float}
\usepackage{color}
\usepackage{url}

\usepackage{xcolor}

\usepackage{siunitx}

\newcommand{\Chalf}{C_{\footnotesize\mbox{half}}}
\newcommand{\Cmem}{C_{\footnotesize\mbox{mem}}}
\newcommand{\Edrug}{E_{\footnotesize\mbox{drug}}}
\newcommand{\Etol}{E_{\footnotesize\mbox{drug}}}
\newcommand{\Edt}{E_{\footnotesize\mbox{drug+tol}}}

\title{Optimal Dosing Schedules for Substances Inducing Tolerance}

\author[a,1]{Golrokh Nouri}
\author[a,1]{Paul Tupper} 

\affil[a]{Department of Mathematics, Simon Fraser University, Burnaby, BC V5A 1S6, Canada }
\affil[1]{Address correspondence to: \texttt{golrokh\_nouri@sfu.ca}, \texttt{pft3@sfu.ca}}
\date{\today}

\begin{document}

\maketitle

\begin{abstract}
Many drugs used therapeutically or recreationally induce tolerance: the effect of the substance decreases with repeated use. This phenomenon may reduce the efficacy of the substance unless dosage is increased beyond what is healthy for the individual. Restoring the effect of the substance can often be obtained by taking a break from consumption. We propose designing dosing schedules that maximize the desired effect of the substance with a given total consumption, while factoring in the effect of tolerance. We provide a simple mathematical model of response to consumption and tolerance that can be fit from data on substance administration and response. Using this model with given parameters, we determine optimal consumption schedules to maximize a given objective. We illustrate with the example of caffeine, where we provide a schedule of consumption for a user who values the effects of caffeine on all days but needs extra alertness on some days of the week.
\end{abstract}

\section{Introduction}

Drug tolerance is the phenomenon where repeated use of a substance leads to diminished effects in an individual. The rate at which tolerance develops depends on the drug, the frequency and dose with which it is administered, the metabolism of the individual, and many other factors. Tolerance is important in medicine, as it means that drug treatments with constant dosage may become less effective over time. Clinicians have to decide when to increase dosage and by how much, balancing the negative side effects of the substance. Likewise, tolerance plays a role in addiction: initially
a substance provides a pleasurable response in a user, 
but over time, the same dose
no longer produces the same effect, and the negative effects of withdrawal occur if consumption is stopped altogether \cite{MolecularBasis,Tolerancedef_merckmanual}.

Drug tolerance can be  acute or chronic, depending on the time-span in which it develops, with acute tolerance occurring after even a single dose (tachyphylaxis), and chronic occurring over longer periods of consumption \cite{MolecularBasis,tolerance_in_anesthesiology_book}.   Tolerance can be categorized in other ways as well. Pharmacokinetic tolerance refers to 
the body developing ways to prevent as much of the substance from reaching the targeted area, such as by clearing it faster from the body \cite{dumas2008opioid}.
Pharmacodynamic tolerance, on the other hand, develops as a result of the body’s homeostatic response to chemical imbalances, for example by down regulating the receptors that pair to a particular substance. Another categorization is based on the different biological levels in play since drug tolerance can affect molecular, cellular, or behavioural processes in the body. Regardless of the specific type or level at which tolerance works, in this paper we focus on the broad features of tolerance at whatever level it occurs, a reduction in response to a drug due to repeated exposure \cite{MolecularBasis,tolerance_in_anesthesiology_book}.

Tolerance is familiar to many people from their experience with the most commonly consumed recreational drug in the world: caffeine, consumed in the form of tea, coffee, and soft drinks \cite{CaffeineConsumptionStats}.  
The effects of caffeine on first consumption are increased blood pressure, improvement in vigilance, and increased alertness and ability to sustain attention, especially in low-arousal situations like early mornings  
\cite{EffectsOfCaffeine_Behaviour,EffectsOfCaffeine_Health}.
It does so by blocking adenosine receptors in the brain \cite{caffmechanismofaction}.
Adenosine is a naturally occurring
compound released by the brain to promote sleep and suppress arousal.
But with a steady daily consumption of coffee the energizing effects of caffeine no longer occur, and in fact when coffee is stopped for a day, the person has the negative side effects of withdrawal such as fatigue and headaches \cite{EffectsOfCaffeine_Health,lara_caffeine}. 
Tolerance occurs for caffeine because the body produces more adenosine receptors, a phenomenon called receptor up-regulation \cite{Ammon1991cafftolerance}. 
Tolerant people therefore need a higher caffeine
concentration to block these receptors, in order to obtain the same level of effect as before.

Psychology has a theory of the development of tolerance, opponent process theory, that can be applied to any stimulus, not just the use of substances
%One theory explains the development of tolerance is the opponent process theory
\cite{solomon1974opponent,solomon1973opponent}. According to this theory, when a subject is exposed to a stimulus that induces an affective reaction, there is  a balancing opponent process initiated that is slower than the initial response.  This opponent process causes the reaction to diminish while the stimulus is still present, and for there to be an opposite reaction after the stimulus is removed.
%"every process that has an affective balance (i.e. is pleasant or unpleasant) is followed by a secondary, "opponent process". This opponent process sets in after the primary process is quieted. 
With repeated exposure, the primary process remains of the same magnitude
%becomes weaker 
while the opponent process becomes stronger. The effect is that after repeated exposure to the stimulus the primary reaction appears weaker and the opposite reaction after the stimulus is removed becomes stronger.
We think of opponent process theory as a higher level theory that encompasses many more detailed theories of tolerance.
Our model of tolerance is inspired by the description in opponent process theory, though we show that it can be applied by fitting to individual data for different drugs.

Our first contribution is to provide a simple mathematical model of reaction to substance that is consistent with the well-known facts of tolerance to substance use. Our model includes two different mechanisms for tolerance. The first mechanism accounts for acute tolerance, and is borrowed from Porchet et al.\ \cite{porchet1988pharmacodynamic}. It accounts for the effect where after a first dose gives an effect of a given size, the second dosage of the same size gives a reduced effect on the user. This mechanism does not lead to a withdrawal or rebound symptoms.  The second mechanism accounts for long-term tolerance, and is inspired by opponent process theory.  This is where continual use of a substance affects the user so that a change is measurable even when the substance has not recently been used.  In particular, this mechanism accounts for withdrawal.
%that is inspired by opponent process theory and consistent with the well-known facts of tolerance to substances. We introduce a model that captures both acute and chronic drug tolerance, depending on the parameter setting.  The model is a system of four ordinary differential equations that together describe effect of a substance on blood concentration and its subsequent effects on tolerance.
Both mechanisms can be turned off or on independently of each other by an appropriate choice of parameters.
Though for most substances we expect both mechanisms to be present, 
 we demonstrate the applicability of the model with only one mechanism at a time.
We use two different parameters settings, one describing acute tolerance to nicotine \cite{porchet1988pharmacodynamic}, and another describing long-term tolerance to caffeine \cite{lara_caffeine}.
 
 Our second contribution shows how once we have a set of parameters for our model, we can use optimization to design optimal dosing schedules for a substance. 
 For example, a subject may have chronic pain, and use an analgesic to provide relief. Continual use without breaks may lead to a state where the subject's experience of pain is not much less than originally without the drug. Does providing breaks in the dosing schedule make sense as a strategy to minimize tolerance?
 We do a detailed exploration with caffeine, imagining a coffee drinker who wants to be alert on particular days of the week more than others, and has a maximum total amount of caffeine that they want to consume per week. 
 
 %We use this model to experiment with different dosing regimes and to optimize treatment plans for a desired schedule. Our model is then fit to data from a nicotine dose-response trial to gather appropriate parameters. \pt{I think we do nicotine and then caffeine separately?}

%Given the existence of tolerance, how can we design dosage plans that optimize a desired outcome? 
 %If the pain fluctuates with time....?
%Our second contribution is to show how to optimize a dosage regime for a particular substance. Our example is a coffee drinker who wants to be alert on particular days.

Our model is similar to those developed by Peper \cite{peper2004theory,peper2004theory2} and Porchet et al. \cite{porchet1988pharmacodynamic}, though in particular it is simpler than Peper's.  The simplicity allows us to determine easily interpretable parameters from available data, and means that the computational demands of simulating the data are light enough to enable optimization of dosing strategies.

\section{Model}
The model is a system of differential equations together describing the response of a subject to an administered substance over time.
%substance tolerance over some time period.

First, we let
$D(t)$ represent the rate at which the drug is administered at time $t$, in units of mass per time. Often, in practice $D(t)$ will be zero most of the time except for some short intervals where it is positive, when the substance is actually being consumed.
We let the function $C(t)$ represent the blood plasma concentration of the substance at time $t$ in units of mass per volume.  %The change in blood plasma concentration 
The dynamics of $C(t)$ over time is 
%at any given time is  
given by:
\begin{equation} \label{eq1}
\frac{d}{dt}C(t) = k_{1}( k_{7} D(t) - C(t) )
\end{equation}
Here, the first term $k_1 k_7 D(t)$ represents the administration of the drug increasing the blood plasma level, while the  $k_{1} C(t)$ term represents the excretion of the substance from the body.
%is subtracted from the drug administered to account for the distribution of the substance in the digestive tract. 
A constant dose $D$ leads to an equilibrium blood concentration of $C = k_7 D$, and
$k_{1}$ determines the rate at which the system approaches this equilibrium.
%of the digestive system’s response to drug administration, and
%$k_{7}$ is its sensitivity to one dose of the drug.

For tolerance to develop over time, the body needs some measure of memory of the substance previously administered. We model this with the quantity $\Cmem(t)$ where:
\begin{equation}
\label{eq:Cmem}
\frac{d}{dt}\Cmem(t) = k_{5} ( C(t) - \Cmem(t) )
\end{equation}
% $\Cmem(t)$ is a measure of the average history of the blood concentration level in the body; 
% In the Supplementary information we show that 
For $\Cmem(0)=0$ it is straightforward to check that  
 \[
 \Cmem(t) = k_5^{-1} \int_0^t e^{-k_5  (t-s)}  C(s) \, ds
 \]
 so $\Cmem(t)$ can be interpreted as an average of $C(t)$ over the time leading up to $t$, with more recent values being more heavily weighted.
% of how much substance is still left in the body, and
The parameter $k_{5}$ controls how long this memory is, with $\Cmem$ being a weighted average of $C(t)$ over a time period of order $1/k_5$.  $\Cmem(t)$ could either represent the concentration of the substance in some longer lasting reservoir than the blood, or some other
consequences of the body’s memory of the substance, for example, the number of extra receptors grown for the substance.
%We can think of $\Cmem(t)$ as the area under the $C(t)$ curve, weighted exponentially heavier as t approaches current time. The derivation for \eqref{eq:Cmem}  can be found in the supplementary materials section.

The quantity $E(t)$ is the effect under consideration of the substance at time $t$. For example, for caffeine, $E(t)$ may be alertness; for an analgesic, $E(t)$ may be a level of pain.  Our model for $E(t)$ as a function of time is 
\begin{equation} \label{eqn:Edynamics}
E(t) = E_b(t) + \frac{F(t)}{ 1+ \Cmem(t) / \Chalf}
\end{equation}
where $F(t)$ is an idealized effect of the drug (without tolerance) and $E_b(t)$ is the current baseline value of $E$, each of which have their own dynamics which we will describe below.

Our two mechanisms of tolerance work to push $E(t)$ away from desired values and can be seen in \eqref{eqn:Edynamics}. 
The first mechanism, for acute tolerance, acts through the denominator of the second term on the right. If $\Cmem=0$ then this term is at its maximum value: $F(t)$. Any larger value of $\Cmem$ reduces this term, and therefore makes $E(t)$ closer to $E_b(t)$. The strength of this effect is controlled by $\Chalf$. When $\Cmem=\Chalf$ then the effect of $F$ is halved. We describe the dynamics of $F$ below. If we want to turn off this mechanism, we set $\Chalf=\infty$.
The second mechanism, for long-term tolerance, works through the term $E_b(t)$. Below we will describe the dynamics of $E_b$, where exposure to the drug shifts its value in order to move $E(t)$ away from more desirable values. 

%For a given subject, we hypothesize that there is an underlying baseline level of this quantity $E_0$. In the long-term absence of the substance, $E(t)$ has $E_0$ as its equilibrium value. Then we assume there is an idealized effect of the drug $F(t)$ and current baseline $E_b(t)$.

%The primary effect of the substance is to move  $E(t)$ away from its baseline level and toward a more desirable value. 
%We model the effects of tolerance in two different ways. 
%We model the time course of $E(t)$ via

%$E_b$ is a baseline value of $E$ that is effected by prior consumption of the drug. $F$ is the idealized underlying effect of the substance
% if there were no tolerance. 

We model the dynamics of $F$ by 
\begin{equation}
\label{eq3}
\frac{d}{dt}F(t) = k_2( k_6 C(t) - F(t) )
\end{equation}
where $k_2$ is the rate at which $F$ converges to its value determined by $C$, and $k_6$ determines the relationship between concentration of the drug $C$ and underlying effect $F$.
 
%The second mechanism, for long-term tolerance, acts through 
%The denominator in the fraction has the effect of decreasing the effect of the drug: when $\Cmem=\Chalf$ the effect of the drug is half as great.  
%But the substance has two effects on $E(t)$: one direct and one indirect. The direct effect of the substance to shift $E(t)$ to a new (more desirable) value, and occurs through $C(t)$, blood plasma concentration of the substance. The indirect effect is that of tolerance, and reduces the effect of the drug by pushing $E(t)$ towards $E_0$, and occcurs through $\Cmem(t)$. 
%We model this through assuming that there is a variable baseline value of $E(t)$, $E_b(t)$, that equilibriates to $E_0$ in the absence of the substance but is affected by $\Cmem(t)$:
We model the dynamics of $E_b$ by
\begin{equation} 
\label{eq4}
\frac{d}{dt}E_{b}(t) = k_{3} ( E_{0} - k_{4} \Cmem(t) - E_{b}(t) ).
\end{equation}
For a given subject, we hypothesize that there is an underlying baseline level of this quantity $E_0$. 
%In the long-term absence of the substance ($\Cmem=0$), $E(t)$ has $E_0$ as its equilibrium value. 
Note that if $\Cmem=0$, so that the drug and any memory of it is cleared from the body, $E_b(t)$ converges to $E_0$ with rate $k_3$. Likewise, if $\Cmem(t)$ is constant, $E_b$ converges to $E_0 - k_4 \Cmem$, so $k_4$ determines just how much memory of the substance shifts the baseline value of $E$.
Setting $k_4$ to zero turns off the long-term tolerance mechanism.

Note that through appropriate parameter choices, our model can incorporate both tolerance mechanisms ($k_4>0$, $\Chalf<\infty$) just long-term tolerance ($k_4>0$, $\Chalf=\infty$), just acute tolerance ($k_4=0$, $\Chalf<\infty$), or neither ($k_4=0$, $\Chalf=\infty$).

\subsection{A single rapid dose, with no tolerance}
First we see what our model does in the situation where a dose is delivered very quickly, and we do not consider tolerance of either sort.
This corresponds to $\Chalf=\infty$ and $E_b = E_0$.
Letting $\tau$ be a short interval of time, and $\mathcal{D}$ be the total dose, we have $D(t) = \mathcal{D}/\tau$ for $t \in [0,\tau]$ and $0$ otherwise. Roughly, neglecting the decay term, $dC/dt = k_1 k_7 D(t)$ over time interval $[0,\tau]$ and so if $C(0)=0$ we get $C(\tau)=\tau ( k_1 k_7 \mathcal{D}/\tau) = k_1 k_7 \mathcal{D}$. After time $\tau$, $D$ is zero, and so the equation reduces to $dC/dt = -k_1 C(t)$. Hence $C(t)= k_1 k_7 \mathcal{D} e^{-k_1 (t-\tau)}$. Letting $\tau$ go to zero gives 
\begin{equation} \label{eq:fastdose1}
C(t) =  k_1 k_7 \mathcal{D} e^{-k_1 t}
\end{equation}
after a dose given at time 0.
With neither form of tolerance included in the model, we have that $E(t) = E_0 + F(t)$ where $F(t)$ satisfies \eqref{eq3}.
%\[
%dE/dt = k_2 ( E_0 + k_6 C(t) - E(t) )
%\]
If $k_2$ is large, then $E$ will rapidly converge to $E_0 + k_6 C(t)$ on short time intervals. If we use the small $t$ approximation in \eqref{eq:fastdose1} we get that $E$ is approximately $E_0 + k_6 k_1 k_7 \mathcal{D}$ at its peak immediately after the dose has taken its effect, and then decays to $E_0$ at approximately the same rate as $C(t)$.
%, where we have used the small $t$ approximation in \eqref{eq:fastdose1}.

\subsection{Constant Dosage, Long-Term Tolerance} 
In order to get rough analytical results for the development of long-term tolerance over time, we consider the case where $D(t)=D$ is constant. %Equivalently, $D$ may very rapidly, but on a time scale that is much shorter than $k_1^{-1}$, so that the effect on $C(t)$ in equation \ref{eq1} is that of an average.
%This may provide insert text about averaging.
We let $\Chalf=\infty$ so that $E=E_b + F$.
Then $C(t)$ converges to $k_7 D$ over a time period of scale $k_1^{-1}$.
In this case $\Cmem(t)$, being an average of $C(t)$, will also converge to $k_7 D$, and assuming $k_5 \ll k_1$, this will occur on the time scale $k_5^{-1}$. 
The same applies for setting $D$ to zero after a period of substance use. Blood plasma concentration will return to $0$ on a time scale of $k_1^{-1}$ and $\Cmem$ will return to $0$ on a time scale of $k_5^{-1}$. 

Now suppose that $\Cmem$ has equilibrated to $k_7 D$ and remains constant. Then $E_b$ will equilibriate to $E_0- k_4 \Cmem = E_0 - k_4 k_7 D$ over a time of scale $k_3^{-1}$. Likewise, $E(t)$ will equilibriate to $E_b + k_6 C = E_0 - k_4 k_7 D + k_6 k_7 D$. This first term $E_0$ is the subject's intrinsic base line, the second term is the effect of tolerance, and the third term is the direct effect of the drug. For most substances we expect $k_4, k_6>0$ so that the drug increases $E$, but the effect of tolerance is to push in the opposite direction. Depending of the relative size of $k_4$ and $k_6$, the net effect of the drug (including tolerance) can increase or decrease $E$.

To a rough approximation we can see that these equilibrium calculations fit the basic phenomenology of substance use and tolerance. Before tolerance has time to develop, steady use of the substance changes the value of $E$ from $E_0$ to $\Edrug = E_0 + k_6 k_7 D$. Over the longer term, tolerance develops and the new value for $E$ is $\Edt = E_0 - k_4 k_7 D + k_6 k_7 D$, which is reduced relative to the initial effect of the drug. Now, if the drug is suddenly removed, $E$ goes to the value $\Etol= E_0 - k_4 k_7 D$, which is even worse than without the drug at all. However, eventually the effect of tolerance will wear off, and $E$
will return to $E_0$.

These considerations give us a formula for estimating $k_4$. 
\begin{equation} \label{eqn:k4_set}
k_4= k_6 \frac{ \Edrug - \Edt}{\Edrug - E_0}
\end{equation}

\section{Examples}
\par 
We now explore the behaviour of our model for two sets of parameter values; one selected to capture the effects of nicotine consumption over the course of a few hours \cite{porchet1988pharmacodynamic}, another that of caffeine over many days \cite{lara_caffeine}.
%We can use our model to see what effects different dosing regimens have on the subject factoring in the effect of tolerance.
%We can now use this model to observe the behaviour of tolerance under different dosing strategies of caffeine consumption.

 \subsection{Nicotine}
\begin{figure}[ht!]
\centering
\includegraphics[width=1\linewidth]{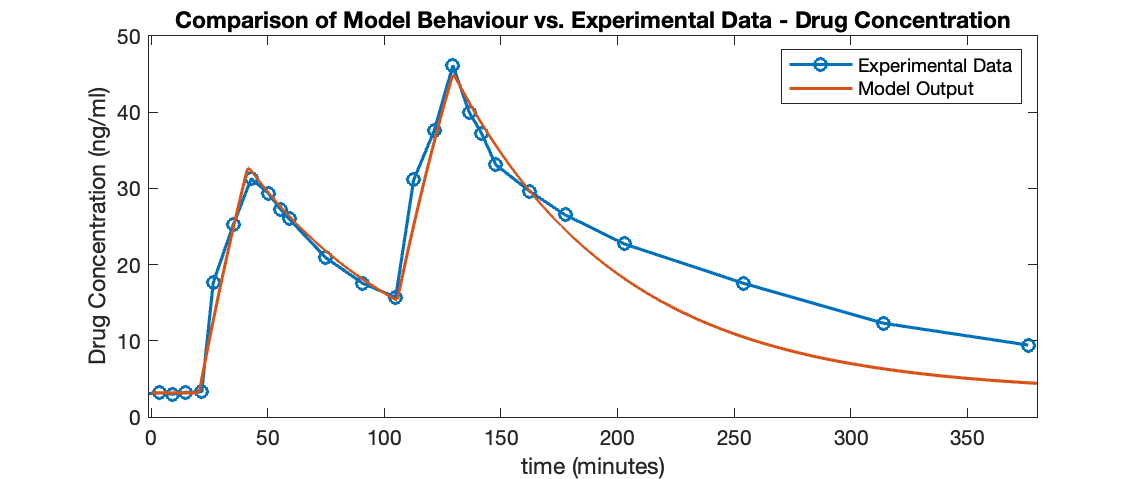}
\includegraphics[width=1\linewidth]{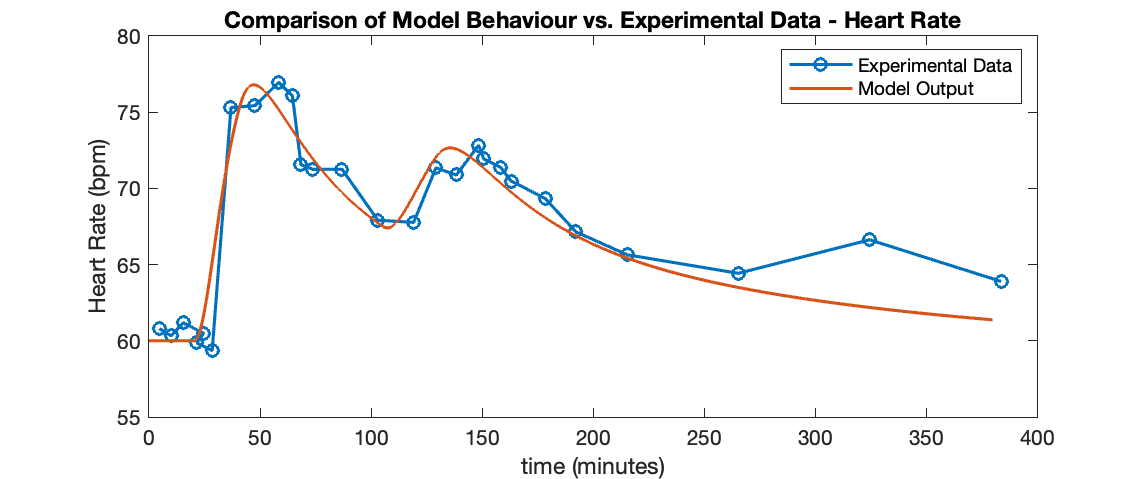}
\caption{A fit of the model to experimental data from Porchet et al. \cite{porchet1988pharmacodynamic}. Parameters are listed in Table~\ref{table: CafNicParameters}. }
\label{fig2}
\end{figure}

First we demonstrate our model with the development of acute (over a short time period) tolerance with nicotine. Porchet et al.\ \cite{porchet1988pharmacodynamic} administered nicotine intravenously with two doses in short succession to a subject, and measured the effect on blood concentration of nicotine and on heart rate. One effect of nicotine was to increase the subjects' heart rate from a baseline value.   They found that the second dose led to a lower peak in heart rate compared to the first, indicating the development of acute tolerance. We attempted to fit this data with our model. 
Since the experiment occurred over a few hours,  we turned off our longer-term tolerance mechanism by setting $k_3=k_4=0$, and hence $E_b(t)=E_0$ for all time. 
We selected other parameters by hand in order to get a reasonable fit to Porchet et al.'s data (the average measurements over eight subjects). The parameters are shown in Table~\ref{table: CafNicParameters} and the fits are shown in Figure~\ref{fig2}.

We see that we are able to capture the main features of the experimental data. Importantly, we see that there is no withdrawal or rebound effect after the consumption of all the nicotine: the heart rate does not drop below the inital heart rate of the subject. Earlier we attempted to capture the short term acute tolerance evident here with a shift in the value of $E_b$. We found that this always led to a withdrawal effect that is not apparent in this data. Nicotine withdrawal does indeed indeed lead to lower than normal heart rate (bradycardia) \cite{dani2011neurophysiology}, but we could not find time series data showing this.
 
% \par In our first example, we use our model to observe the development of acute tolerance to nicotine administration. 
% \par We estimate the parameters corresponding to nicotine from a set of experimental data found in a paper by Prochet et. al. \cite{porchet1988pharmacodynamic}.  
 
 %In this paper, a different model of drug tolerance is introduced and tested on data measuring the relationship between consecutive doses of nicotine and heart rate. In this experiment, 8 healthy individuals are given 2 i.v. infusions of nicotine, administered one hour apart, and their heart rate is recorded. Since the data from the experiment was not available, 
 %we used MATLAB's GRABIT tool to 
 %we extracted the data points from the figures in the paper. We then by hand fitted our model to this data in order to test its overall performance, as well as establish estimates for the model's parameters.  The estimated parameters from the fitting are found in the last row of table (\ref{table: CafNicParameters}). Figure \ref{fig2} also illustrates the comparison between nicotine concentrations and heart rate in the experiment, and $C(t)$ and $E(t)$ plotted using our model and these parameters.

 The main shortcoming of our model in fitting $C(t)$ as a function of time is that our model predicts a faster decay of concentration after both doses are complete. Porchet et al. are able to capture this better with their similar model, as they includes a second compartment in the pharmacokinetic model, something that could be added to our model if desired. 
 
% \pt{Paul will add something about our poorer fit, the lack of a second compartment, and some qualitative features. Also note that this is acute tolerance, and probably there is a separate chronic tolerance mechanism for nicotine.}

\begin{table}[ht]
\centering
\caption{Parameters used in caffeine and nicotine example}
\begin{tabular}{lrrrrrrrrr}
Parameter & $E_{0}$ & $k_{1}$  & $k_{2}$ & $k_{3}$ & $k_{4}$ & $k_{5}$ & $k_{6}$ & $k_{7}$ & $\Chalf$ \\ 
\midrule
Units &    & 1/min & 1/min & 1/day & mL/$\mu$g & 1/day & mL/$\mu$g  & min/ mL & $\mu$g/mL \\
\midrule
Nicotine values & 60 & 0.014 & 0.08 & 0 & 0 & 20. & 1.8 $\times 10^3$ & 0.0175 & 0.005 \\
\midrule
Caffeine values & 0  & 0.002  & 0.1 & 0.5   & 0.3   &  0.5  &  0.4  &  0.0125 &  $\infty$ \\
\bottomrule
\label{table: CafNicParameters}
\end{tabular}
\end{table}

\subsection{Caffeine}
We next consider longer term effects of caffeine consumption.  We chose parameters %$E_{0}, k_{1},...,k_{7}$ 
to roughly correspond to the effects of caffeine on a typical person. For example, peak plasma caffeine concentration from consuming $\mathcal{D} = $100 \si{mg} of caffeine (equivalent to about 1 cup of coffee \cite{FDAcaffeine}) has a value of about $C= $ 2.5 \si{\mu g/mL} and then the concentration decays with a half-life of about 5 hours,
%happens about 30 minutes after consumption, 
according to a study by Teekachunhatean et al.\ \cite{PlasmaConc1}. This means that the rate of decay is approximately $k_1 = \ln 2 / (300$\si{min}) =0.002 \si{min^{-1}}.
If we then select  $k_7$ in Equation~\ref{eq:fastdose1} to match these values we get  $k_7 = $2.5 \si{\mu g/mL} $/(k_1  10^{5}$ \si{\mu g}$)=$0.0125 \si{min/mL}.  
Caffeine has a variety of effects, and rather than selecting one we take $E(t)$ to be some form of ``alertness'',  a dimensionless quantity going up to a peak value of near 1, and having a baseline of 0. This sets the value of $k_6= 1/( 2.5$ \si{\mu g/ mL})$=0.4$ \si{mL/\mu g}. We do not observe big differences reported between times of peak concentration of caffeine and peak psychological effect \cite{corti2002coffee}, so we assume that it is relatively fast, and set $k_2=$0.1 \si{min^{-1}}. 

To estimate parameter values related to the development of tolerance for caffeine, we look to Lara et al.\ \cite{lara_caffeine}. 
They study the effect of caffeine on exercise, and how this diminishes after repeated use. We summarize the results of one of their experiments as shown in \cite[Fig.\ 2]{lara_caffeine}. Looking at this figure we obtain rough estimates of $\Edrug$ and $\Edt$ as follows. 
After a caffeine free period, athletes' performance in an exercise test was measured; in researchers' normalized units $E_0=0$. Then they began a daily routine of caffeine consumption followed by the same exercise test.  Initially there is a significant improvement in performance $\Edrug \approx 2$, but this declines over time, about halving over 10 days. We estimate that the eventual performance under continuing caffeine would be $\Edt \approx 0.5$. Then after 10 days, the daily caffeine is stopped, and the subjects' performance on that day declines to be less than it was initially without caffeine (appoximately -0.5), indicating a withdrawal effect. 
%Tolerance develops on a time scale of 10 or so days. 

Since their measurements are only once per day, we don't have any way to measure acute tolerance. Hence we set $\Chalf=\infty$, therefore turning of that mechanism.
The data doesn't allow us to distinguish between the time constants for $\Cmem$ and $E_b$,
so we assume $k_3=k_5= 1/ (2 \, $\si{days}$)= 0.5$ \si{days^{-1}}.  %We choose $E_0=0$.
Using our estimates for $E_0, \Edrug, \Edt$ above in  \eqref{eqn:k4_set} we get  that $k_4= k_6 1.5/ 2= 0.3 $\si{mL/\mu g}.
The complete list of the parameters can be found in Table \ref{table: CafNicParameters}. 

We use these parameters to explore the behaviour of our model for caffeine.  We introduce three different dosing regimens and observe the changes in blood caffeine concentration $C(t)$, alertness $E(t)$, and the changes in the baseline alertness $E_b(t)$ for each regime. 
%Let $E(t)$ be a measure of ``alertness and energy levels'' as a percentage, and let it have a baseline value of 60. 
For all three dosing regimens, the experiment will run for a total of 4 weeks, and the total amount of caffeine consumed will be fixed at 2800 mg (28 cups of coffee) over the whole duration of the experiment.

\par
\begin{figure*}[ht]
\centering
\includegraphics[width=0.32\linewidth]{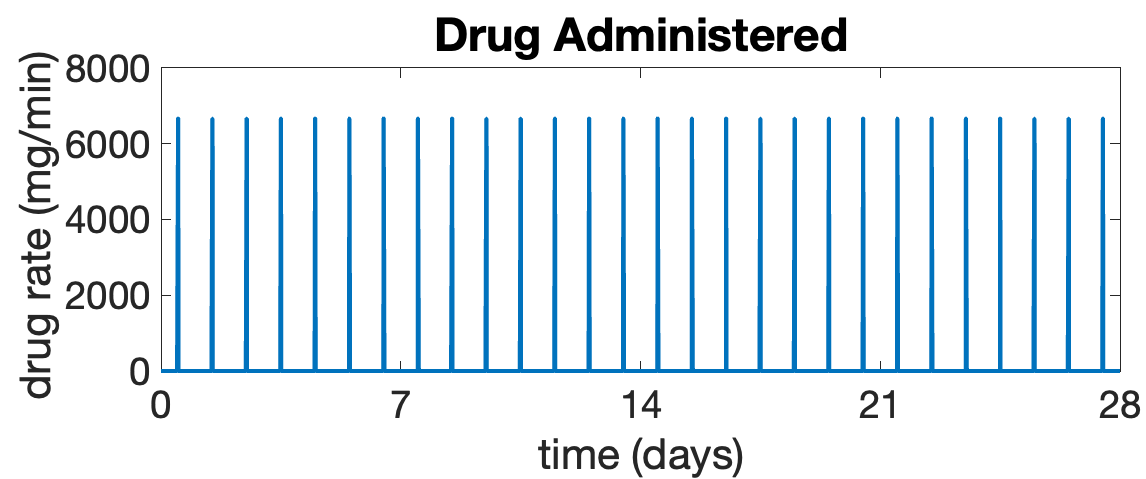}
\includegraphics[width=0.32\linewidth]{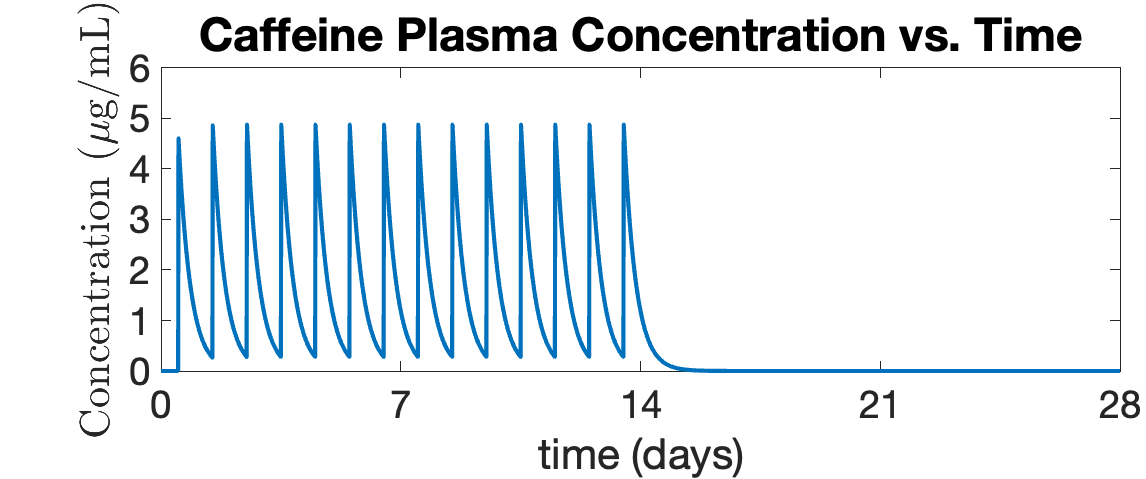}
\includegraphics[width=0.32\linewidth]{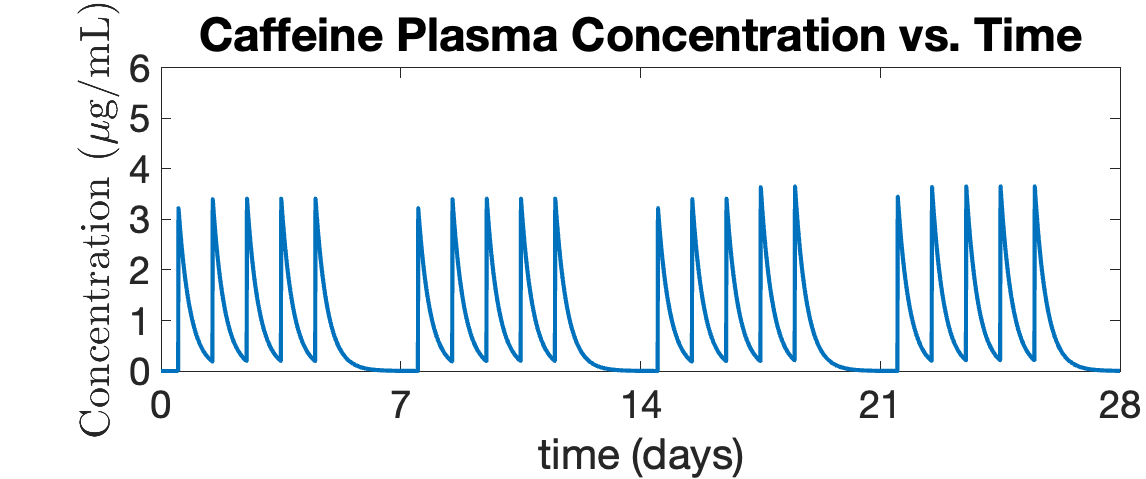}
\includegraphics[width=0.32\linewidth]{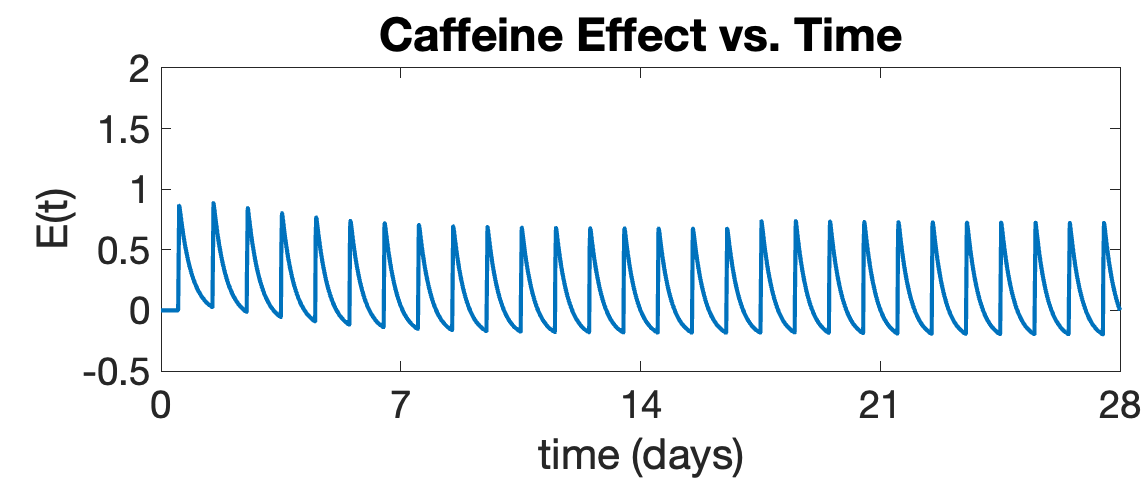}
\includegraphics[width=0.32\linewidth]{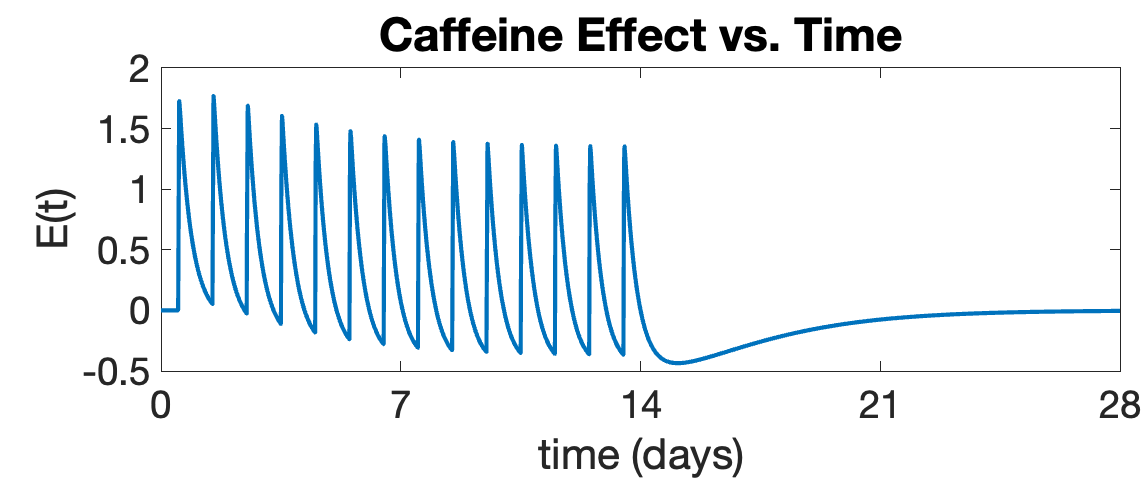}
\includegraphics[width=0.32\linewidth]{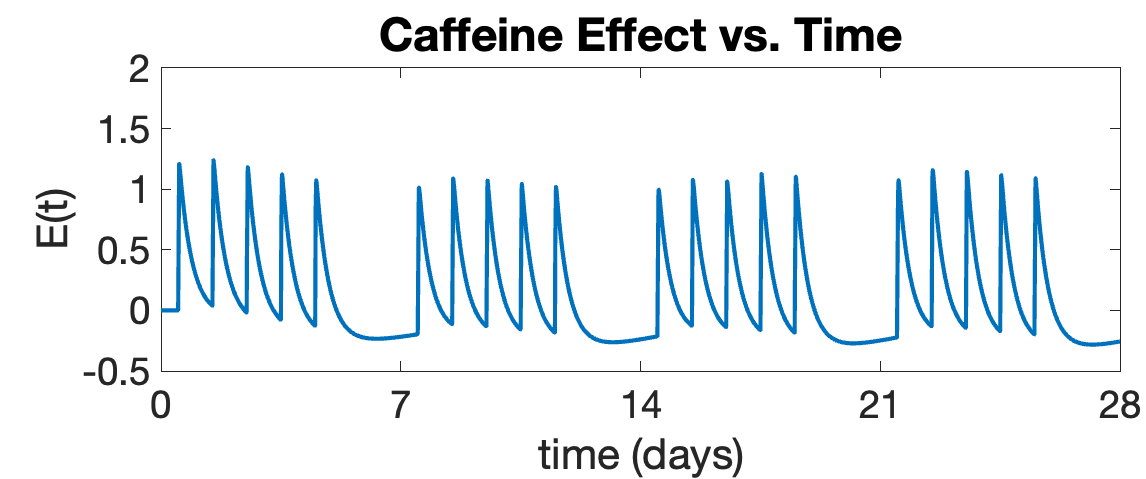}
\includegraphics[width=0.32\linewidth]{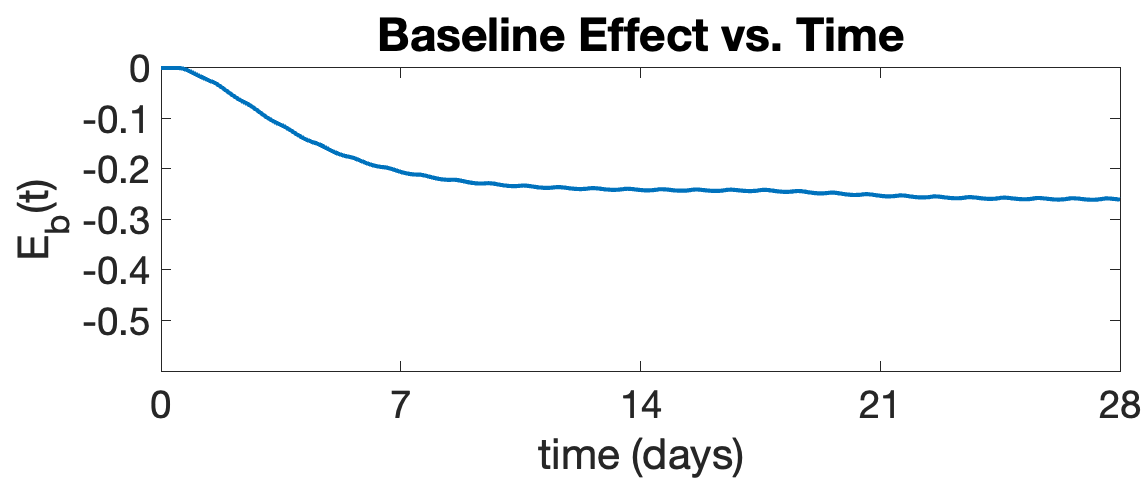}
\includegraphics[width=0.32\linewidth]{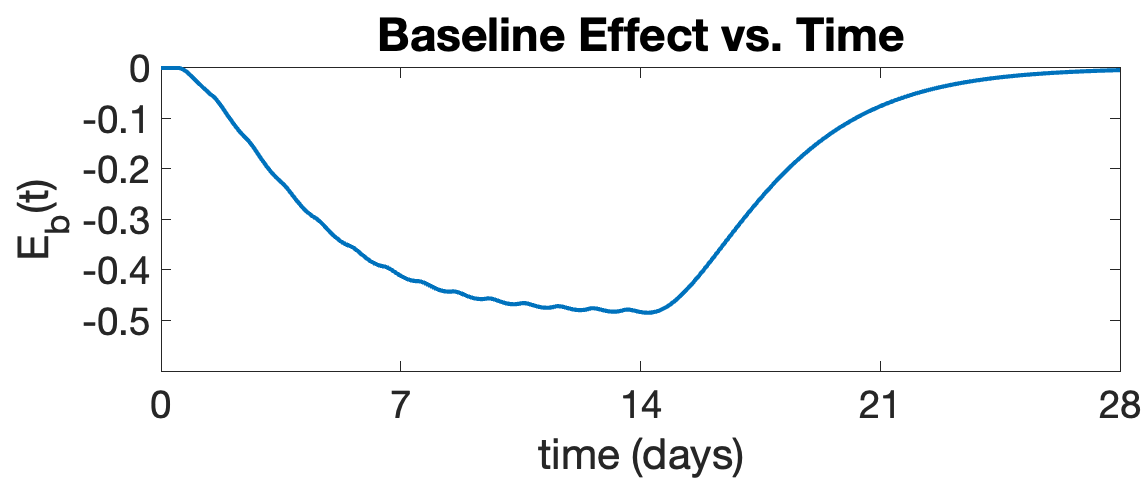}
\includegraphics[width=0.32\linewidth]{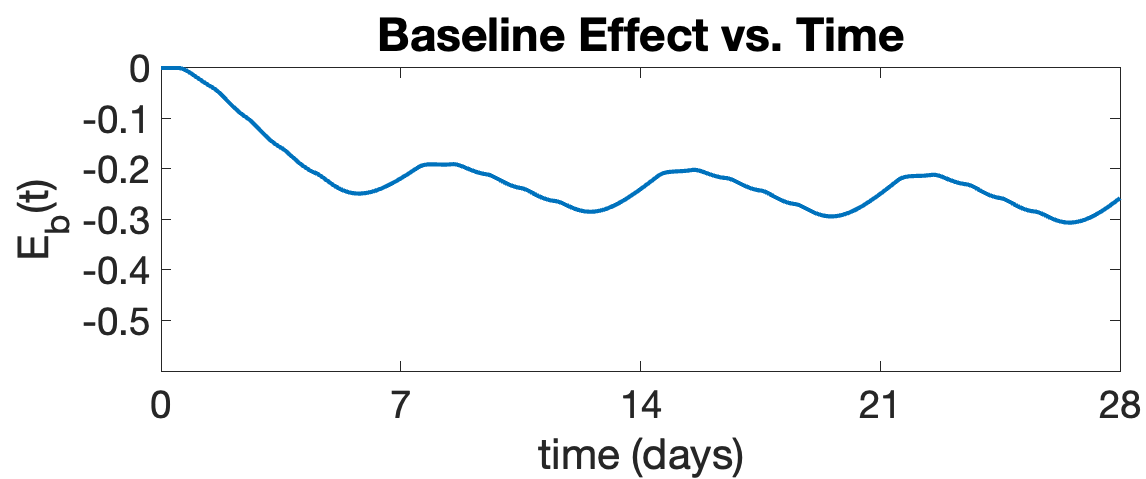}
\caption{Caffeine plasma level $C(t)$, alertness $E(t)$, and baseline alertness $E_b(t)$ over a four-week period for three dosing regimens.}
%This will be helpful to compare the parameters.}
\label{fig1}
\end{figure*}

 The first dosing regimen (Figure \ref{fig1}, left) is to consume one cup of coffee (100mg) on a daily basis for the entirety of the experiment. The subject consumes it starting at noon over a 15 minute period. We observe that the effect $E(t)$ peaks shortly after consumption each day. The highest peak occurs on the first day of consumption. Peak height lowers on each subsequent day, assymptoting to a lower value after a couple of weeks. This is reflected in $E_b$ decaying from 0 to approximately  $-0.25$ and remaining there. Another feature is that $E$ is eventually negative at times during the day.

\par In the second regime (Figure \ref{fig1}, centre) the subject has two cups of coffee per day for two weeks and then stops completely. The higher dose leads to higher peaks, which similarly lessens in height with each day. When the subject stops consumption, alertness drops to below 0 consistently for many days. But this slowly wears off and alertness returns to the zero baseline.

%\par Next we will look at how taking long tolerance breaks can affect the baseline levels. The second dosing regimen will be to drink 2 cups of coffee every day for 2 weeks, followed by a 2 week hiatus. 

\par The last dosing regime we will look at is to consume 140 $mg$ of caffeine on weekdays, followed by no consumption over the weekend in order to prevent tolerance developing. This does indeed lead to a higher peak on the Monday of the second week. This is partly due to there being a higher consumption on that day in this regime, but also the diminished effect of tolerance from consumption on the weekend.

%\par Figure \ref{fig1} shows changes in plasma caffeine concentration, alertness levels, and homeostatic baseline of alertness levels under the three regimens described above. In Figure 1.a, corresponding to the first regimen, note the continuous drop in $E(t)$ levels until day 10, where both the overall alertness level, and the homeostatic baseline level, $E_b(t)$, plateau and a new tolerance to this particular dosing regimen of caffeine is established. The second dosing regimen (figure 1.b) shows that after a new tolerance has been established, the homeostatic response of the body will try to revert back to the original baseline levels, given a long enough tolerance break.

\section{Optimization of Dosing Schedule}
We next consider choosing the optimal dosing regime to obtain a desired outcome. 
As in our previous example shown in Fig.\ \ref{fig1}, we imagine a subject who consumes caffeine in the form of coffee, with the goal of being more alert at specific times during the week.
We use our model with the parameters for caffeine we determined earlier, as shown in Table~\ref{table: CafNicParameters}.
%Let coffee be the substance of interest, and alertness its desired effect. 
Let us assume that our subject has important meetings on Mondays and Thursdays and thus needs to be very alert on those days, but enjoys the effect of caffeine on all days of the week. As in the simulations shown in Figure~\ref{fig1}, the subject is given caffeine every day from noon to 12:15pm. We select the strength of the dose $d_i$ on day $i$ for $i=1,\ldots,7$, in units where 1 indicates a single cup of coffee. Starting from a state with no caffeine in the system $C=\Cmem=E_b=E=0$ we run the simulation for three weeks. 
%Each day of the week has its own dose which is repeated every week. 
%We choose the dose for each day of the week and simulate using our model for three weeks.
We choose the dose on each day to maximize the following objective function. 
On the last week of the simulation, on each day at 3pm we measure the alertness $E$ of the subject, $e_i$ for $i=1,\ldots,7$. We take the weighted average of the alertness measurements over the week
\begin{equation}
     f =   \sum_{i=1}^{7} {w_{i} \sqrt{e_{i}}} 
\end{equation} 
We model our subject's desire to be more alert on Mondays and Thursdays by setting $w_1=w_4=10$  and all other days of the week $w_i=0.2$.
We use the square root in the objective function to model diminishing marginal benefit from increasing levels of alertness; otherwise, the optimum schedule is always to consume as much coffee as possible on the day with the greatest weight assigned to it.

If there are no constraints on the dose given every day of the week, then $f$ can be made arbitrarily large by giving the subject more and more coffee. So we consider two different constraints on coffee consumption. In the first, the \emph{daily constraint}, the subject doesn't consume more than 2 cups of coffee a day. In the second, the \emph{weekly constraint}, the subject doesn't consume for than 10 cups of coffee each week.

%If we asked them to rate the importance of alertness from the scale of 0 to 10 on those days, they'd rate Mondays and Thursdays a 10, and other days of the week a 1. The subject also doesn't want to consume more than 10 cups of coffee each week, due to health reasons.
%We aim to use optimization tools on our model of drug tolerance to produce the ideal weekly schedule of caffeine consumption whose outcome is similar to the subject's needs.
%To maximize $f$ given these constraints, we used Matlab's \textit{fmincon} applied to $-f$. The input is a vector of length 7, $d_i$, $i=1,\ldots,7$ giving the dose per day of the week in units of cups of coffee.
Initially, we used an adaptive solver for our system of differential equations, which makes use of changing step sizes to produce more accurate results more efficiently. This worked well when testing the model and fitting it to data, but not when used in optimization, since the adaptive step sizes led to the computed $f$ not being a continuous function of the dosing schedule. 
%Solvers with adaptive step sizes can produce discontinuous outputs in the form of rises and dips. However, these local minima can be mistaken as the global minimum by the optimizer function. Our goal is to find the dosing schedule that results in the smallest value of the objective function, i.e. the global minimum of the objective function. To find this global minimum we need a smooth function, one that a solver with adaptive step sizes cannot produce \textcolor{violet}{(is this always true?)}. 
So, we used the fixed-step forward Euler method  with a small step size to solve the system of differential equations within the optimization. 

\begin{figure}[ht!]
\centering
\includegraphics[width=1\linewidth]{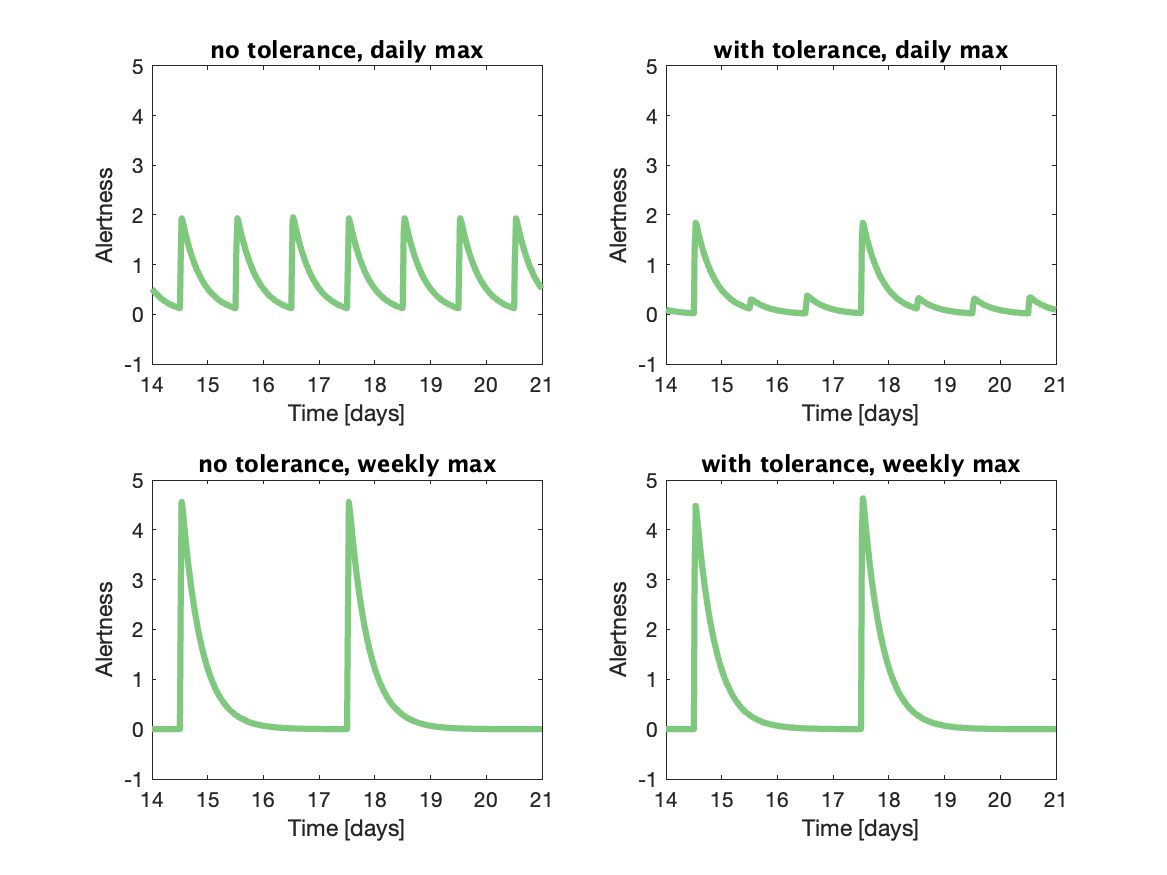}
\caption{Alertness versus time over the course of a week of caffeine consumption with optimized schedules. Left column, results for the model without tolerance; right column, with tolerance. Top row, consumption restricted to at most two cups of coffee a day; bottom row, consumption restrictred to ten cups of coffee a week.  }
\label{fig:optimized}
\end{figure}

We initialized the optimizer with dose $d_i$ for $i=1,\ldots,7$ drawn independently and uniformly at random in $[0,1]$. The optimizer was run to convergence. In Table~\ref{tab:optimal_sched} we show the optimal schedules determined for four different conditions. We varied whether we included long-term tolerance in the model ($k_4=0.3$ mL/$\mu$g) or not ($k_4=0$), and whether we used the daily constraint or the weekly constraint. In Table~\ref{tab:optimal_sched} we show the optimal schedule for the four different conditions. In Figure~\ref{fig:optimized} we show the resulting alertness $E$ versus time in the last week of the simulation with the optimized schedule. 

The results without tolerance are  straightforward to understand. With a daily maximum it is always best to consume as much coffee on every day. This naturally yields the same alertness peak on each day. With a weekly max it is best to split almost all coffee consumption equally between the two days when it is needed most, though interestingly 1/1000 of a cup of coffee is recommended on Sundays, presumably because some alertness is desired on the other days, and Sunday is the day with the least residual effect from Monday and Thursday.

Results change when we add tolerance. With a daily maximum the optimal choice is to maximize coffee consumption to 2 cups on Monday and Thursday, but to have significantly reduced consumption on other days of the week:  between 20\% and 40\% of a cup on these days. The model determines this to be the optimal tradeoff between being alert on these days and not having tolerance deprive the subject of the benefits of caffeine on the important days. With the weekly maximum now all the consumption occurs on Monday and Thursday. Slightly more caffeine is recommended on Thursday than Monday, and no caffeine consumption at all is recommended on other days.

\begin{table}[ht]
\centering
\caption{Optimal dosing schedule for four different conditions}
\begin{tabular}{lrrrrrrr}
Condition & Day 1 & Day 2 & Day 3 & Day 4 & Day 5 & Day 6 & Day 7 \\ 
\midrule
no tolerance, daily max &   2.0000 &   2.0000  &  2.0000  &  2.0000  &  2.0000  &  2.0000   & 2.0000 \\
\midrule
with tolerance, daily max &   2.0000 &  0.2179  &  0.3921 &   2.0000  &  0.2456 &   0.3275  &  0.3568 \\
\midrule
no tolerance, weekly max &  4.9997  &  0.0000 &   0.0000 &   4.9990  &  0.0000  &  0.0000  &  0.0013 \\
\midrule
with tolerance, weekly max & 4.9180 &   0.0000  &  0.0000 &   5.0820   &  0.0000 &   0.0000  &  0.0000 \\
\bottomrule
\label{tab:optimal_sched}
\end{tabular}
\end{table}

% [3.156300753649019,0.936199509130598,0.555925629080220,3.324416280343012,1.058130786836284,0.644361282231676,0.324665737909160]

\section{Discussion}

We have presented a simple model of the effect of a drug on an individual and the subsequent development of tolerance. We fit the model to two sets of data, one showing the effect of acute tolerance to nicotine and the other showing long-term tolerance to caffeine. Our model includes two different mechanisms for tolerance, one for acute and one for long-term. We found that we could not fit both datasets with a model with a single mechanism. In particular, our acute tolerance mechanism does not exhibit withdrawal, and so cannot match what was observed for caffeine \cite{lara_caffeine}. On the other hand, our long-term tolerance model based must exhibit withdrawal after removal of the substance, and so cannot match the data for nicotine from Porchet et al \cite{porchet1988pharmacodynamic}. We expect that most substances will exhibit both of these forms of tolerance, and direction for future work is to derive parameters for the model for a substance incorporating both of these effects, once such data is available.

As with any model there are phenomena that it will not be able to account for. Most obviously from our fit to nicotine data above, our single compartment pharmacokinetic model prevents us from accurately capturing the decay of nicotine concentration over time. This could be fixed easily (at the cost of adding more parameters). A more serious problem is the difficulty of getting sufficient data to fit all the parameters. Above we used a combination of fits to data and educated guesses. But in reality our model, and our choices of parameters, need validation against a wide range of experimental data before being used as a basis for application in physiology. 

Our second contribution is to show that given such a model with appropriate parameters we can use optimization to determine optimal dosing schedules for particular goals.  We only considered a case where long-term tolerance (over the course of days) was relevant, but our methods are in principle applicable to more short term situations where acute tolerance is more important. For example, if a subject is competing in a day-long chess tournament, determining what is the optimal consumption of nicotine or caffeine for performance.

We chose the familiar example of caffeine to illustrate our method for dosing schedule optimization, but there are many more important instances that come to mind, such as in medicine and in addiction management. 
Ideally such considerations could inform treatment design for these problems. But the soundness of the recommendations relies on the soundness of the underlying model. There is also the problem of turning what is wanted from a dosing regime (such as freedom from pain) into an objective function to be optimized.
Our hope though, is that by using models such as ours, and exploring optimal dosing schedules in different regimes, ideally in collaboration with a clinician, we will be able to learn something of the issues involved in managing drug treatment.

\section*{Data and code availability}
All data and code for this study are available in the GitHub repository:\\ {\tt{https://github.com/PaulFredTupper/optimal-dosing-schedules}}

\section*{Acknowledgements}
This study was funded by the Natural Science and Engineering Research Council (Canada) Discovery Grant (RGPIN-2019-06911).

% Bibliography
\bibliographystyle{plain}

\end{document}